\documentstyle[12pt,epsf]{ioplppt}

\begin{document}

\jl{1}

\title{$q$-Deformed Fock spaces and modular representations of spin symmetric groups}[
$q$-Deformed Fock spaces]

\author{Bernard Leclerc\dag\   and Jean-Yves Thibon\ddag}

\address{\dag D\'epartement de Math\'ematiques, Universit\'e de Caen,
Esplanade de la Paix, BP 5186, 14032 Caen cedex, France}

\address{\ddag Institut Gaspard Monge, Universit\'e de Marne-la-Vall\'ee,
2 rue de la Butte-Verte, 93166 Noisy-le-Grand cedex, France}


\begin{abstract}
We use  the Fock space representation of the quantum
affine algebra of type $A^{(2)}_{2n}$ to obtain a description
of the global crystal basis of its basic level 1 module.
We formulate a conjecture relating this basis to decomposition
matrices of spin symmetric groups in characteristic $2n+1$.
\end{abstract}


\def\SG{{\rm S}}
\def\SP{{\rm \widehat{S}}}
\def\C{{\bf C}}
\def\Z{{\bf Z}}
\def\Q{{\bf Q}}
\def\slchap{\widehat{sl}}
\def\Des{{\rm Des\,}}
\def\maj{{\rm maj\,}}
\def\mod{{\ \rm mod\ }}
\def\<{\langle}
\def\>{\rangle}
\def\pr#1{\langle #1 \rangle}
\def\prh#1{\pr{\widehat{#1}}}
\def\vac{|0\>}
\def\DP{{\rm DP}}
\def\DPR{{\rm DPR}}
\def\SR{{\rm SR}}
\def\ket#1{|#1\>}

\def\shuff#1#2{\mathbin{
      \hbox{\vbox{
        \hbox{\vrule
              \hskip#2
              \vrule height#1 width 0pt
               }%
        \hrule}%
             \vbox{
        \hbox{\vrule
              \hskip#2
              \vrule height#1 width 0pt
               \vrule }%
        \hrule}%
}}}

\def\shuffle{{\mathchoice{\shuff{7pt}{3.5pt}}%
                        {\shuff{6pt}{3pt}}%
                        {\shuff{4pt}{2pt}}%
                        {\shuff{3pt}{1.5pt}}}}%
\def\shuf{\,\shuffle\,}
\def\qshuf{\shuf_q}

\newtheorem{example}{Example}[section]
\newtheorem{theorem}[example]{Theorem}
\newtheorem{corollary}[example]{Corollary}
\newtheorem{definition}[example]{Definition}
\newtheorem{proposition}[example]{Proposition}
\newtheorem{algorithm}[example]{Algorithm}
\newtheorem{lemma}[example]{Lemma}
\newtheorem{conjecture}[example]{Conjecture}

\font\FontSeteight=msbm8
\font\FontSetnine=msbm9
\font\FontSetten=msbm10
\font\FontSettwelve=msbm10 scaled 1200
\newdimen\Squaresize \Squaresize=18pt
\newdimen\Thickness \Thickness=0.5pt

\def\Square#1{\hbox{\vrule width \Thickness
   \vbox to \Squaresize{\hrule height \Thickness\vss
      \hbox to \Squaresize{\hss#1\hss}
   \vss\hrule height\Thickness}
\unskip\vrule width \Thickness}
\kern-\Thickness}

\def\Vsquare#1{\vbox{\Square{$#1$}}\kern-\Thickness}
\def\blk{\omit\hskip\Squaresize}

\def\young#1{
\vbox{\smallskip\offinterlineskip
\halign{&\Vsquare{##}\cr #1}}}

\def\thisbox#1{\kern-.09ex\fbox{#1}}
\def\downbox#1{\lower1.200em\hbox{#1}}

\section{Introduction}

The Fock space representation of the
quantum affine algebra $U_q(\slchap_n)=U_q(A^{(1)}_{n-1})$ 
was constructed by Hayashi \cite{H}.
A combinatorial version of this construction was then used by
Misra and Miwa \cite{MM} to describe  Kashiwara's crystal basis of
the basic representation $V(\Lambda_0)$.
This made it possible to compute the global crystal basis of
$V(\Lambda_0)$ \cite{LLT}. Then, it was conjectured
that the degree $m$ part of the transition matrices 
giving the coefficients of the global basis 
on the natural basis of the Fock space
were $q$-analogues of the decomposition matrices of the  type $A$
Hecke algebras $H_m$ at an $n$th root of unity \cite{LLT}. 
According to a conjecture
of James \cite{J}, these should coincide, for $n$ prime and
large enough,
with the decomposition matrices of symmetric groups $\SG_m$ over a field
of characteristic $n$.
The conjecture of \cite{LLT} has been proved by Ariki \cite{Ar},
and by Grojnowski \cite{Gr} using the results of \cite{G}.

There is another approach to the calculation of decomposition
matrices of type $A$ Hecke algebras, relying upon Soergel's
results on tilting modules for quantum groups at roots of
unity \cite{Soe1,Soe2}.
This approach also leads to $q$-analogues of decomposition
numbers expressed in terms of Kazhdan-Lusztig polynomials.
It seems that these $q$-analogues are the same as those
of \cite{LLT} but there is no proof of this coincidence.
In fact, the relationship between the two approaches is still unclear.

The results of \cite{LLT,Ar,Gr} have been applied recently
by Foda {\it et al.} \cite{FLOTW} to determine which simple
$H_m$-modules remain simple after restriction to $H_{m-1}$
and to show that this problem is equivalent to the decomposition
of a tensor product of level 1 $A_{n-1}^{(1)}$-modules.
This provided an explanation for an intriguing correspondence
previously observed in \cite{FOW} between a class of RSOS models
and modular representations of symmetric groups.

Another description of the $U_q(A^{(1)}_{n-1})$-Fock space,
as a deformation of the infinite wedge realization of
the fermionic Fock space, was obtained by Stern \cite{St}.
In \cite{KMS}, the $q$-bosons needed for the decomposition
of the Fock space into irreducible $U_q(A^{(1)}_{n-1})$-modules
were introduced. This construction was used in \cite{LLTrib}
to give a combinatorial formula for the highest weight
vectors, and in \cite{LT} to define a  canonical basis
of the whole Fock space which was conjectured to
yield the decomposition matrices
of $q$-Schur algebras at roots of unity.
Moreover, strong support in favor of this conjecture was
obtained by establishing its compatibility with a version 
of the Steinberg tensor product theorem proved by James
in this context \cite{J,LT}.

Recently, the theory of perfect crystals \cite{KMN1,KMN2} allowed
Kashiwara {\it et al.} \cite{KMPY} to define a general
notion of $q$-Fock space, extending the results of \cite{KMS}
to several series of affine algebras.
Their results apply in particular to the twisted affine algebra
of type $A^{(2)}_{2n}$, which is the case considered in this note.

It has been noticed by Nakajima and Yamada \cite{NY} that the combinatorics
of the basic representation
$V(\Lambda_n)$ of $A^{(2)}_{2n}$ was similar to the
one encountered in the $(2n+1)$-modular representation theory of the spin
symmetric groups $\SP_m$ by Morris \cite{Mo1} as early as 1965.
This can be explained to a certain extent by observing that
the $(r,\bar{r})$-inducing operators of Morris and Yaseen \cite{MY}
coincide with the Chevalley lowering operators of the
Fock space representation of $A^{(2)}_{2n}$. This provides
a further example of the phenomenon observed in \cite{LLT}
in the case of symmetric groups and $A_{n-1}^{(1)}$-algebras. 

In this note, we give the analogues for $U_q(A^{(2)}_{2n})$ of the
results of \cite{LLT}. 
Using the level~1 $q$-Fock spaces of \cite{KMPY},
we describe an algorithm for computing
the canonical basis of the basic representation $V(\Lambda_n)$, which
allows us to prove that this basis is in the $\Z[q]$-lattice
spanned by the natural basis of the $q$-Fock space, and that
the transition matrices have an upper triangle of zeros
(Theorem 4.1). 

We conjecture that the specialization $q=1$ gives, up to splitting of rows and
columns for pairs of associate characters, and for sufficiently
large primes $p=2n+1$, the decomposition matrices of spin symmetric groups. 
However, the reduction $q=1$ is more tricky than in the $A_{n-1}^{(1)}$ case.
Indeed, the $q$-Fock space of $A^{(2)}_{2n}$ is strictly larger than
the classical one, and one has to factor out the null space
of a certain quadratic form \cite{KMPY} to recover the usual
description.

The missing ingredient in the spin case when we compare it to
\cite{LLT} is that, since the spin symmetric groups are
not Coxeter groups, there is no standard way of associating
to them a Hecke algebra, and this is an important obstruction
for proving our conjecture.
What we can actually prove is that all self-associate 
projective characters of $\SP_m$ are linear combinations
of characters obtained from smaller groups by a sequence
of $(r,\overline r)$-inductions (Theorem 6.1). 
This proof is constructive in the sense that the intermediate
basis $\{A(\mu)\}$ of our algorithm for the canonical basis,
suitably specialized at $q=1$, is a basis for the space spanned 
by such characters.

This should have implications on the way of labelling the
irreducible modular spin representations of $\SP_m$.
Up to now, a coherent labelling scheme has been found
only for $p=3$ \cite{BMO} and $p=5$ \cite{ABO}.
The case $p\ge 7$ led to formidable difficulties.
To overcome this problem, we propose to use the labels
of the crystal graph of $V(\Lambda_n)$, which may contain
partitions with repeated parts not arising in the
representation theory of $\SP_m$, and corresponding to  ghost vectors
of the $q$-Fock space at $q=1$.

\section{The Fock space representation of $U_q(A^{(2)}_{2n})$}

The Fock space representation of the affine Lie algebra
$A^{(2)}_{2n}$ can be constructed by means of its
embedding in $b_\infty=\widehat{go}_\infty$, the completed infinite
rank affine Lie algebra of type $B$  \cite{DJKM1,DJKM2}.

The (bosonic) Fock space of type $B$ is
the polynomial algebra ${\cal F} = \C[p_{2j+1}, j\ge 0 ]$ in an infinite
number of  generators $p_{2j+1}$ of odd degree $2j+1$. If one identifies
$p_k$ with the power sum symmetric function $p_k=\sum_i x_i^k$
in some infinite set of variables, the natural basis of weight
vectors for $b_\infty$ is given by Schur's $P$-functions $P_\lambda$
(where $\lambda$ runs over the set $\DP$ of partitions
into distinct parts) \cite{DJKM1,You,JY}.

The Chevalley generators $e^\infty_i$, $f^\infty_i$ ($i\ge 0$)
of $b_\infty$ act on $P_\lambda$ by
\begin{equation}\label{FP}
e^\infty_i P_\lambda = P_\mu \ , \qquad f^\infty_i P_\lambda = P_\nu
\end{equation}
where $\mu$ (resp. $\nu$) is obtained from $\lambda$ by replacing its part $i+1$
by $i$ (resp. its part $i$ by $i+1$), the result being $0$  
if $i+1$ (resp. $i$) is not a part of $\lambda$.
Also, it is
understood that $P_\mu=0$ as soon as $\mu$ has a multiple part.
For example, $f^\infty_0 P_{32}=P_{321}$,
$f^\infty_3 P_{32}=P_{42}$,
$e^\infty_1 P_{32}= P_{31}$ and $e^\infty_2 P_{32}=P_{22}=0$.

Let $h=2n+1$.
The Chevalley generators $e_i$, $f_i$ of $A^{(2)}_{2n}$
will be realized as
\begin{equation}\label{FNP}
f_i=\sum_{j\equiv n\pm i} f^\infty_j  \qquad
(i=0,\ldots ,n) \,,
\end{equation}
\begin{equation}
e_i=\sum_{j\equiv n\pm i } e^\infty_j \qquad 
(i=0,\ldots, n-1)\,,\qquad 
e_n=e^\infty_0 + 
2 \sum_{\scriptstyle j>0 \atop \scriptstyle j \equiv 0,-1} e^\infty_j \,,
\end{equation}
where all congruences are taken modulo $h$.
Let $A_{2n}^{(2)}{}'$ be the derived algebra of $A_{2n}^{(2)}$
(obtained by omitting the degree operator $d$).
The action of $A_{2n}^{(2)}{}'$ on ${\cal F}$ is centralized
by the Heisenberg algebra generated by the operators
$\displaystyle{\partial\over\partial p_{hs}}$ and $p_{hs}$ for odd $s\ge 1$.
This implies that the Fock space decomposes under $A_{2n}^{(2)}$ as
\begin{equation}\label{DEC1}
{\cal F} = \bigoplus_{k\ge 0} V(\Lambda_n-k\delta)^{\oplus p^* (k)}
\end{equation}
where $p^*(k)$ is the number of partitions of $k$ into odd parts.
In particular, the subrepresentation generated by the vacuum vector
$|0\rangle=P_0 = 1$ is the basic
representation $V(\Lambda_n)$ of $A_{2n}^{(2)}$, and its principally
specialized character is \cite{KKLW}
\begin{equation}\label{CHAR}
{\rm ch}_t\,V(\Lambda_n) =
\sum_{m\ge 0}\dim V(\Lambda_n)_m\,t^m =
\prod_{\scriptstyle i \ {\rm odd} \atop \scriptstyle i\not\equiv 0\mod h}
{1\over 1-t^i}\,.
\end{equation}

The $q$-deformation of this situation has been discovered
by Kashiwara {\it et al.}  \cite{KMPY}. 
Contrary to the case of
$A^{(1)}_{n-1}$, the $q$-Fock space is strictly larger than
the classical one. 
We recall here briefly their construction, referring to
\cite{KMPY} for details and notation.

Let $\DP_h(m)$ be the set 
of partitions $\lambda=(1^{m_1}2^{m_2}\ldots r^{m_r})$
of $m$ for which $m_i\le 1$ when $i\not\equiv 0 \mod h$. 
For example, $\DP_3(7)=\{(7),(61),(52),(43),(421),(331)\}$.
Set $\DP_h=\bigcup_m \DP_h(m)$.
Then, the $q$-Fock space of type $A_{2n}^{(2)}$ is
\begin{equation}
{\cal F}_q = \bigoplus_{\lambda\in \DP_h} \Q(q)\, |\lambda\>
\end{equation}
where for $\lambda=(\lambda_1,\ldots,\lambda_r)$,
$|\lambda\>$ denotes the infinite $q$-wedge product
\[
|\lambda\> = u_\lambda = u_{\lambda_1}\wedge_q u_{\lambda_2}\wedge_q\cdots\wedge_q
u_{\lambda_r}\wedge_q u_0 \wedge_q u_0 \wedge_q \cdots 
\]
of basis vectors $u_i$ of the representation $V_{\rm aff}$.
The quantum affine algebra $U_q(A_{2n}^{(2)})$ acts on
$V_{\rm aff}=\bigoplus_{i\in\Z}\Q(q) u_i$ 
by 
\begin{eqnarray}
f_i u_j = \cases{ u_{j+1} & if $j\equiv n\pm i\mod h$ \\
                      0   & otherwise \\} 
\qquad (i=0,\ldots,n-1) 
\\
\label{ACTF}
f_n u_j = \cases{ u_{j+1} & if $j\equiv -1 \mod h$ \\ 
       (q+q^{-1}) u_{j+1} & if $j\equiv 0 \mod h$ \\
                        0 & otherwise \\} 
\\
e_i u_j = \cases{ u_{j-1} & if $j\equiv n+1\pm i\mod h$ \\
                     0    & otherwise \\} 
\qquad (i=0,\ldots,n-1)
\\
e_n u_j = \cases{ u_{j-1} & if $j\equiv 1 \mod h$ \\
       (q+q^{-1}) u_{j-1} & if $j\equiv 0 \mod h$ \\
                        0 & otherwise \\}
\\
t_0 u_j = \cases{ q^4 u_j & if $j\equiv n\mod h$ \\
              q^{-4}  u_j & if $j\equiv n+1 \mod h$ \\
                      u_j & otherwise \\}  
\\
t_i u_j = \cases{ q^2 u_j & if $j\equiv n\pm i\mod h$ \\
              q^{-2}  u_j & if $j\equiv n+1\pm i \mod h$ \\
                      u_j & otherwise  \\ }
\qquad (i=1,\ldots,n-1)
\\
t_n u_j = \cases{ q^2 u_j & if $j\equiv -1\mod h$ \\
              q^{-2}  u_j & if $j\equiv 1 \mod h$ \\
                      u_j & otherwise  \\} 
\end{eqnarray}
The only commutation rules we will need to describe the
action of $e_i$ and $f_i$ on ${\cal F}_q$ are:
\begin{eqnarray}
u_j \wedge_q u_j &=& 0 \ {\rm if}\ j\not\equiv 0 \mod h \\
u_j \wedge_q u_{j+1} &=& -q^2 u_{j+1}\wedge_q u_j \ {\rm if} \label{STR2}
 j\equiv 0,-1 \mod h \ .
\end{eqnarray}
The action on the vacuum vector 
$\vac = u_0\wedge_qu_0\wedge_q\cdots $
is given by
\begin{equation}
e_i\vac = 0, \qquad
f_i\vac = \delta_{i n}|1\>, \qquad
t_i\vac = q^{\delta_{i n}}\vac,
\end{equation} 
and on a $q$-wedge 
$|\lambda\>=u_{\lambda_1}\wedge_q\cdots\wedge_q u_{\lambda_r}
\wedge_q \vac$,
\begin{eqnarray}
f_i |\lambda\>
=&
f_iu_{\lambda_1}\wedge_q t_iu_{\lambda_2}\wedge_q\cdots t_iu_{\lambda_r}
\wedge_q t_i\vac \nonumber \\
& +
u_{\lambda_1}\wedge_q f_iu_{\lambda_2}\wedge_q\cdots t_iu_{\lambda_r}
\wedge_q t_i\vac  \nonumber \\
& + \cdots +
u_{\lambda_1}\wedge_q u_{\lambda_2}\wedge_q\cdots u_{\lambda_r}
\wedge_q f_i\vac
\end{eqnarray}
\begin{eqnarray}
e_i |\lambda\>
=&
t_i^{-1} u_{\lambda_1}\wedge_q t_i^{-1}u_{\lambda_2}\wedge_q\cdots t_i^{-1}u_{\lambda_r}
\wedge_q e_i\vac \nonumber \\
& +
t_i^{-1}u_{\lambda_1}\wedge_q t_i^{-1}u_{\lambda_2}\wedge_q\cdots e_iu_{\lambda_r}
\wedge_q \vac  \nonumber \\
& + \cdots +
e_i u_{\lambda_1}\wedge_q u_{\lambda_2}\wedge_q\cdots u_{\lambda_r}
\wedge_q \vac
\end{eqnarray}
\begin{equation}
t_i |\lambda\> =
t_i u_{\lambda_1}\wedge_q t_i u_{\lambda_2}\wedge_q\cdots\wedge_q
t_i u_{\lambda_r}\wedge_q t_i\vac \ .
\end{equation}
For example, with $n=2$, one has
\[
f_2 |542\>= (q^4+q^2)|642\>+q|552\>+|5421\>,
\]
and 
\[
f_2 |552\> = (q^2+1)(|652\>+|562\>)+|5521\> 
= (1-q^4)|652\>+|5521\>,
\]
the last equality resulting from (\ref{STR2}).

It is proved in \cite{KMPY} that ${\cal F}_q$ is an integrable
highest weight $U_q(A_{2n}^{(2)})$-module whose decomposition
into irreducible components, obtained by means of $q$-bosons, is 
\begin{equation}\label{DEC2}
{\cal F}_q = \bigoplus_{k\ge 0} V(\Lambda_n-k\delta)^{\oplus p(k)}
\end{equation}
where $p(k)$ is now the number of all partitions of $k$ 
(compare (\ref{DEC1})).
Thus, the submodule $U_q(A_{2n}^{(2)}) \,|0\>$ is a realization
of the basic representation $V(\Lambda_n)$.

\section{The crystal graph of the $q$-Fock space}

The first step in computing the global basis
of $V(\Lambda_n) \subset {\cal F}_q$ is to determine
the crystal basis of ${\cal F}_q$ whose description
follows from \cite{KMPY,KMN1,KMN2}.
Let $A$ denote the subring of $\Q(q)$ consisting of 
rational functions without pole at $q=0$.
The crystal lattice of ${\cal F}_q$ is 
$L = \bigoplus_{\lambda\in \DP_h} A\,|\lambda\>$,
and the crystal basis of the $\Q$-vector space $L/qL$ is 
$B=\{|\lambda\> \mod qL, \lambda \in \DP_h\}$. 
We shall write $\lambda$ instead of $|\lambda\> \mod qL$.

The Kashiwara operators $\tilde{f}_i$ act on $B$ in 
a simple way recorded on the crystal graph $\Gamma({\cal F}_q)$. 
To describe this graph, one starts with the crystal graph
$\Gamma(V_{\rm aff})$ of $V_{\rm aff}$. This is the graph with vertices 
$j\in \Z$, whose arrows labelled by $i\in \{0,1,\ldots ,n\}$
are given, for $i \not = n$, by
\[
j \stackrel{i}{\longrightarrow} j+1 \quad \Longleftrightarrow \quad
j \equiv n \pm i \mod h \,,
\]
and for $i=n$ by
\[
j \stackrel{n}{\longrightarrow} j+1 \quad \Longleftrightarrow \quad
j \equiv -1,0 \mod h \,.
\]
Thus for $n=2$ this graph is
\[
\cdots \stackrel{1}{\longrightarrow} -1
\stackrel{2}{\longrightarrow} 0
\stackrel{2}{\longrightarrow} 1
\stackrel{1}{\longrightarrow} 2
\stackrel{0}{\longrightarrow} 3
\stackrel{1}{\longrightarrow} 4
\stackrel{2}{\longrightarrow} 5
\stackrel{2}{\longrightarrow} 6
\stackrel{1}{\longrightarrow} 7
\stackrel{0}{\longrightarrow}
\cdots
\]
The graph $\Gamma({\cal F}_q)$ is obtained inductively 
from $\Gamma(V_{\rm aff})$ using the following rules.
Let $\lambda = (\lambda_1,\ldots ,\lambda_r)\in B$,
and write $\lambda = (\lambda_1,\lambda^*)$
where $\lambda^* = (\lambda_2,\ldots ,\lambda_r)$.
Then one has $\tilde{f}_i (0) = \delta_{in} (1)$, 
$\varphi_i(0)= \delta_{in}$, and
\[
\tilde{f}_i\lambda = \left\{
\matrix{(\tilde{f}_i \lambda_1, \lambda^*)
\ {\rm if} \ \varepsilon_i(\lambda_1) \ge \varphi_i(\lambda^*), \cr
(\lambda_1,\tilde{f}_i\lambda^*)
\ {\rm if} \ \varepsilon_i(\lambda_1) < \varphi_i(\lambda^*). } \right.
\]
Here, $\varepsilon_i(\lambda_1)$ means the distance in $\Gamma(V_{\rm aff})$ 
from $\lambda_1$ to the origin of its
$i$-string, and $\varphi_i(\lambda^*)$ means the distance in
$\Gamma({\cal F}_q)$ from $\lambda^*$ to the end of its $i$-string.

Thus for $n=1$ one computes successively the following $1$-strings
of $\Gamma({\cal F}_q)$
\[
(0) \stackrel{1}{\longrightarrow} (1)
\]
\[
(2)=(2,0)\stackrel{1}{\longrightarrow}(2,1) 
\stackrel{1}{\longrightarrow} (3,1) 
\stackrel{1}{\longrightarrow} (4,1) 
\]
\[
(3,2)=(3,2,0)\stackrel{1}{\longrightarrow}(3,2,1) 
\stackrel{1}{\longrightarrow} (3,3,1) 
\stackrel{1}{\longrightarrow} (4,3,1) 
\]
from which one deduces that $\tilde{f_1}(3,3,1) = (4,3,1)$
and $\varphi_1(3,3,1)=1$.
 
\begin{figure}[t]
\begin{center}
\leavevmode
\epsfxsize = 15cm
\epsffile{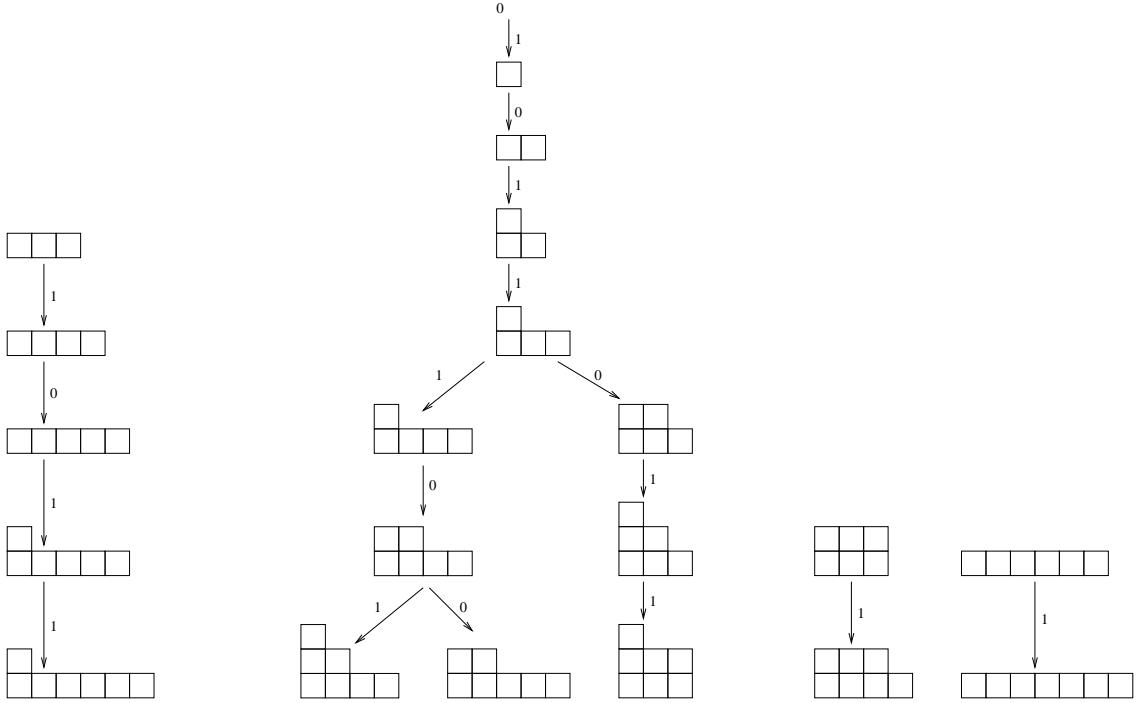}
\end{center}
\caption{\label{FIG1} The graph $\Gamma({\cal F}_q)$ for $A_2^{(2)}$ up to degree $7$} 
\end{figure}

The first layers of the crystal $\Gamma({\cal F}_q)$ for $n=1$
are shown in Fig.~\ref{FIG1}.
One can observe that the decomposition of $\Gamma({\cal F}_q)$ into connected
components reflects the decomposition (\ref{DEC2})
of ${\cal F}_q$ into simple modules. 
More precisely, the connected components of 
$\Gamma({\cal F}_q)$ are all isomorphic as colored graphs
to the component $\Gamma(\Lambda_n)$ containing the
empty partition.
Their highest vertices are the partitions $\nu$ whose
parts are all divisible by $h$. 
This follows from the fact, easily deduced from the
rules we have just explained, that if
$\nu = h\mu = (h\mu_1,\ldots ,h\mu_r)$ is such a partition,
then the map 
\begin{equation}\label{MAP}
\lambda \mapsto \lambda + \nu = (\lambda_1+h\mu_1,\lambda_2+h\mu_2,\ldots \ )
\end{equation}
is a bijection from $\Gamma(\Lambda_n)$ onto the connected component
of $\Gamma({\cal F}_q)$ containing $\nu$, and this bijection commutes with
the operators $\tilde{e}_i$ and $\tilde{f}_i$.
This implies that the vertices of $\Gamma(\Lambda_n)$
are the partitions $\lambda=(\lambda_1,\ldots ,\lambda_r,0)\in \DP_h$ such that
for $i=1,2,\ldots ,r$, one has $\lambda_i- \lambda_{i+1} \le h$ and
$\lambda_i- \lambda_{i+1} < h$ if $\lambda_i \equiv 0 \mod h$.
We shall call a partition that satisfies these conditions
$h$-regular.
The set of $h$-regular partitions of $m$ will be denoted by $\DPR_h(m)$,
and we shall write $\DPR_h=\bigcup_m \DPR_h(m)$.

For example,
\[
\DPR_3(10) = \{ (3331), (4321), (532), (541) \} \,.
\]

\section{The canonical basis of $V(\Lambda_n)$}\label{SECT4}

In this section, we describe an algorithm for computing
the canonical basis (global lower crystal basis) of the
basic representation $V(\Lambda_n)=U_q(A_{2n}^{(2)}) \vac$
in terms of the natural basis $|\lambda\>$ of the
$q$-Fock space. To characterize the canonical basis, we
need the following notations
\begin{equation}
q_i =
\cases{q & if $i=n$ \\
      q^2& if $1\le i<n$ \\
      q^4& if $i=0$ \\ }
\qquad
t_i =
\cases{q^{h_n}& if $i=n$\\
      q^{2h_i}& if $1\le i<n$ \\
      q^{4h_0}& if $i=0$ \\}
\end{equation}
and
\begin{equation}
[k]_i = {q_i^k-q_i^{-k}\over q_i-q_i^{-1}}\ ,
\qquad
[k]_i! = [k]_i [k-1]_i \cdots [1]_i \ .
\end{equation}
The $q$-divided powers of the Chevalley generators are defined by
\begin{equation}
e_i^{(k)} = {e_i^k\over [k]_i!}\ ,\qquad
f_i^{(k)} = {f_i^k\over [k]_i!}\ .
\end{equation}
The canonical basis is defined in terms of an involution
$v\mapsto\overline{v}$ of $V(\Lambda_n)$.
Let $x\mapsto \overline{x}$ be the ring automorphism of
$U_q(A_{2n}^{(2)})$ such that $\overline{q}=q^{-1}$,
$\overline{q^h}=q^{-h}$ for $h$ in the Cartan subalgebra
of $A_{2n}^{(2)}$, and $\overline{e_i}=e_i$,
$\overline{f_i}=f_i$. Then, for $v=x\vac\in V(\Lambda_n)$,
define $\overline{v}=\overline{x}\vac$.

We denote by $U_\Q^-$ the sub-$\Q[q,q^{-1}]$-algebra
of $U_q(A_{2n}^{(2)})$ generated by the $f_i^{(k)}$
and set $V_\Q(\Lambda_n)=U_\Q^-\vac$.
Then, as shown by Kashiwara \cite{K}, there exists a unique
$\Q[q,q^{-1}]$-basis $\{G(\mu), \mu\in \DPR_h\}$ 
of $V_\Q(\Lambda_n)$, such that
\begin{quote}
(G1) $G(\mu) \equiv |\mu\> \mod qL$  

(G2) $\overline{G(\mu)}= G(\mu)$.     
\end{quote}

To compute $G(\mu)$, we follow the same strategy as in
\cite{LLT}. We first introduce an auxiliary basis
$A(\mu)$ satisfying (G2), from which we  manage to construct
combinations satisfying also (G1). More precisely,
let ${\cal F}_q^m$ be the subspace of ${\cal F}_q$ spanned
by $|\lambda\>$ for $\lambda \in \DP_h(m)$ and set 
$V(\Lambda_n)_m={\cal F}_q^m\cap V(\Lambda_n)$. Denote
by $\unlhd$ the natural order on partitions.
Then, the auxiliary
basis will satisfy
\begin{quote}
(A0) $\{A(\mu),\mu\in \DPR_h(m)\}$ is a $\Q[q,q^{-1}]$-basis
of $V_\Q(\Lambda_n)_m$,

(A1) $A(\mu)=\sum_\lambda a_{\lambda\mu}(q)|\lambda\>$,
where $a_{\lambda\mu}(q)=0$ unless $\lambda\unrhd\mu$,
$a_{\mu\mu}(q)=1$ and $a_{\lambda\mu}(q)\in\Z[q,q^{-1}]$,

(A2) $\overline{A(\mu)}=A(\mu)$.
\end{quote}
The basis $A(\mu)$ is obtained
by applying monomials in the $f_i^{(k)}$ to the highest weight vector,
that is, $A(\mu)$ is of the form
\begin{equation}\label{defA}
A(\mu) = f_{r_s}^{(k_s)}f_{r_{s-1}}^{(k_{s-1})}\cdots f_{r_1}^{(k_1)}\vac
\end{equation}
so that (A2) is satisfied. 

The two sequences $(r_1,\ldots,r_s)$ and $(k_1,\ldots,k_s)$ are, as
in \cite{LLT}, obtained by peeling off the $A_{2n}^{(2)}$-ladders
of the partition $\mu$, which are defined as follows. We first fill
the cells of the Young diagram $Y$ of $\mu$ with integers 
(called residues), constant in
each column of $Y$. If $j\equiv n\pm i \mod h$
($0\le i\le n$), the numbers filling
the $j$-th column of $Y$ will be equal to $i$. A ladder of $\mu$
is then a sequence of cells with the same residue, located
in consecutive rows at horizontal distance $h$, except when the residue 
is $n$, in which case two consecutive $n$-cells  in a row belong also
to the same ladder. For example, with $n=3$ and $\mu=(11,7,7,4)$,
one finds $22$ ladders (indicated by subscripts), the longest
one being the 7th, containing three 3-cells:
\[
\young{
3_{19} & 2_{20} & 1_{21} & 0_{22} \cr
3_{13} & 2_{14} & 1_{15} & 0_{16} & 1_{17} & 2_{18} & 3_{19} \cr
3_{7}&2_8&1_9&0_{10}&1_{11}&2_{12}&3_{13}\cr
3_1&2_2&1_3&0_4&1_5&2_6&3_7&3_7&2_8&1_9&0_{10}\cr}
\]
Note that this definition of ladders agrees with that of \cite{BMO}
for $n=1$, but differs from that of \cite{ABO} for $n=2$.

Then, in (\ref{defA}), $s$ is the number of ladders,
$r_i$ the residue of the $i$th ladder, and $k_i$ the number of its
cells. Thus, proceeding with our example,
\[
\fl
A(11,7,7,4)=
f_0f_1f_2f_3^{(2)}f_2f_1f_0f_1f_2f_3^{(2)}f_2f_1f_0^{(2)}
f_1^{(2)}f_2^{(2)}f_3^{(3)} f_2f_1f_0f_1f_2f_3 \vac \ .
\]
The proof of (A0) and (A1) can be readily adapted from
\cite{LLT}. In particular, (A1) follows from the fact that
a partition $\lambda$ belongs to $\DPR_h$ if and only if all cells of a given
ladder intersecting $\lambda$ occupy the highest
possible positions on this ladder.

Another choice of an intermediate basis,
more efficient for practical computations,
would be to use inductively the vectors $G(\nu)$ already computed
and to set $A(\mu)=f_{r_s}^{(k_s)} G(\nu)$, where
$\nu$ is the partition obtained from $\mu$ by removing
its outer ladder. 

Define now the coefficients $b_{\nu\mu}(q)$ by
\begin{equation}
G(\mu) =\sum_\nu b_{\nu\mu}(q) A(\nu) \ .
\end{equation}
Still following \cite{LLT}, one can check that $b_{\nu\mu}(q)=0$
unless $\nu\ge \mu$, where $\ge$ denote the   lexicographic
ordering on partitions, and that $b_{\mu\mu}(q)=1$. Therefore,
one can apply the triangular process of \cite{LLT} as follows. 

Let $\mu^{(1)} < \mu^{(2)} <\ldots < \mu^{(t)}$ be the set
$\DPR_h(m)$ sorted in  lexicographic order, so that
$A(\mu^{(t)})=G(\mu^{(t)})$.  Suppose that the expansion
on the basis $|\lambda\>$ of $G(\mu^{(i+1)}),\ldots, G(\mu^{(t)})$
has already been calculated. Then,
\begin{equation}
G(\mu^{(i)})=
A(\mu^{(i)})-\gamma_{i+1} (q) G(\mu^{(i+1)}) - \cdots - \gamma_t(q)G(\mu^{(t)}) \ ,
\end{equation}
where the coefficients are  determined by the conditions
\[
\gamma_s(q^{-1})=\gamma_s(q), \qquad G(\mu^{(i)}) \equiv |\mu^{(i)}\> \mod qL.
\]

Thus, for $n=1$, the first partition for which $A(\mu)\not = G(\mu)$
is $\mu=(3321)$ and
\begin{eqnarray}
\fl A(3321) =
|3321\> + q|333\> + (q^2-q^6)|432\> + (1+2q^2)|531\> + (q^2+q^4)|54\> 
\nonumber \\
 \lo+ (2q^2+q^4)|621\> +2q^3 |63\> + (q^4+q^6)|72\> +q^4|81\> +q^5|9\>
\end{eqnarray}
Indeed, $A(3321)\equiv |3321\>+|531\> \mod qL$. 
On the other hand, $ A(531)=|531\>+q^2|54\>+q^2|621\>+q^3|63\>+q^6|72\> $
is equal to $G(531)$,
and one finds by subtracting this from $A(3321)$ that
\begin{eqnarray}
\fl G(3321) =  |3321\> + q|333\> + (q^2-q^6)|432\> + 2q^2|531\> + q^4|54\> \\
\lo+ (q^2+q^4)|621\> +q^3 |63\> + q^4|72\> +q^4|81\> +q^5|9\> \ .
\end{eqnarray}
Since $A(432)=|432\>+q^4|531\>+q^2|72\>+q^6|81\>$ satisfies
(G1) and (G2), it has to be equal to
$G(432)$, which completes the determination of the canonical
basis for $m=9$. For $m=10$, the results are displayed as the
columns of Table \ref{TAB1}. 

\begin{table}
\caption{\label{TAB1}
The canonical basis for $n=1$ and $m=10$.
}
\begin{indented}
\item[]\begin{tabular}{@{}llllllll}
\br
&$(3 3 3 1)$&$(4 3 2 1)$&$(5 3 2)$&$(5 4 1)$\\
\mr
$(3 3 3 1)$&1&0&0&0\\
$(4 3 2 1)$&$q-q^{5}$&1&0&0\\
$(4 3 3)$&$q^{2}$&$q$&0&0\\
$(5 3 2)$ &0&0&1&0\\
$(5 4 1)$&$q+q^{3}$&$q^{2}+q^{4}$&0&1\\
$(6 3 1)$&$2\ q^{2}$&$q^{3}$&0&$q$ \\
$(6 4)$&$q^{4}$&0&0&$q^{3}$\\
$(7 2 1)$&$q^{3}+q^{5}$&$q^{2}$&0&$q^{4}$\\
$(7 3)$&$q^{4 }$&$q^{3}$&0&$q^{5}$\\
$(8 2)$ & 0 & 0 & $q^2$ & 0 \\
$(9 1)$&$q^{4}$&$q^{5}$&0&0\\
$(10)$&$q^{6}$&0&0&0\\
\br
\end{tabular}
\end{indented}
\end{table}

In the Fock space representation of $A_{n-1}^{(1)}$, 
the weight of a basis vector $|\lambda\>$ is determined by
the $n$-core of the partition $\lambda$ (and its degree) \cite{ANY,LLT}. 
There is a similar result of Nakajima and Yamada \cite{NY}
for $A_{2n}^{(2)}$, in terms of the notion of $\overline{h}$-core
of a strict partition
introduced by Morris \cite{Mo1} in the context of the modular
representation theory of spin symmetric groups.  

One way to see this is to use a theorem of \cite{MY1}
according to which $\lambda, \mu \in \DP(m)$
have the same $\overline{h}$-core if and only if 
they have, for each $i$, the same number $n_i$ of nodes of residue $i$.
On the other hand, it follows from the 
implementation of the Chevalley generators
that $|\lambda\>$ has $A_{2n}^{(2)}$-weight
$\Lambda_n - \sum_{0\le i \le n} n_i \alpha_i$,
and the statement follows.

The definition of $\overline{h}$-cores can be extended to $\DP_h$
by deciding that if $\lambda$ has repeated parts, its $\overline{h}$-core
is equal to that of the partition obtained by removing those repeated parts.
Then it is clear that if $|\lambda\>$ and $|\mu\>$ have the same 
$U_q(A_{2n}^{(2)})$-weight, the two partitions $\lambda$ and $\mu$ 
have the same $\overline{h}$-core. 
It follows, since $G(\mu)$ is obviously a weight vector, that its
expansion on the basis $|\lambda\>$ involves only partitions
$\lambda$ with the same $\overline{h}$-core as $\mu$.

Summarizing the discussion, we have:

\begin{theorem}\label{TH}
For $\mu \in \DPR_h(m)$, define $d_{\lambda\mu}(q)$ by 
$\displaystyle G(\mu)=\sum_{\lambda \in \DP_h(m)} d_{\lambda\mu}(q)|\lambda\>$.
Then, 

{\rm (i)} $d_{\lambda\mu}(q)\in \Z[q]$,

{\rm (ii)} $d_{\lambda\mu}(q)=0$ unless $\lambda\unrhd\mu$,
and $d_{\mu\mu}(q) = 1$,

{\rm (iii)} $d_{\lambda\mu}(q)=0$ unless $\lambda$ and $\mu$
have the same $\overline{h}$-core.
\end{theorem}

\section{The reduction $q=1$}
As observed by Kashiwara {\it et al.} \cite{KMPY}, to recover the
classical Fock space representation ${\cal F}$ of $A_{2n}^{(2)}$, one has 
to introduce the inner product on ${\cal F}_q$ for 
which the vectors $|\lambda\>$ are orthogonal and the adjoint 
operators of the Chevalley generators are
\begin{equation} \label{ADJOINT}
f_i^{\dag} = q_i e_i t_i, \qquad
e_i^{\dag} = q_i f_i t_i^{-1}, \qquad
t_i^{\dag} = t_i.
\end{equation}
It can be checked that, for $\lambda \in \DP_h$, 
\begin{equation} \label{NORM}
\<\lambda|\lambda\> = \prod_{k>0}\prod_{i=1}^{m_{kh}} (1-(-q^2)^i),
\end{equation}
where $m_{kh}$ is the multiplicity of the part $kh$ in $\lambda$.

Let ${\cal F}_1$ denote the $A_{2n}^{(2)}$-module obtained by specializing
$q$ to 1 as in \cite{KMPY}. This space is strictly larger than the classical Fock space
${\cal F}$, since the dimension of its $m$th homogeneous component
(in the principal gradation) is $|\DP_h(m)|$ whereas that of ${\cal F}$
is only $|\DP(m)|$.
Let ${\cal N} = {\cal F}_1^\perp$ denote the nullspace. It follows
from (\ref{ADJOINT}) that ${\cal N}$ is a $A_{2n}^{(2)}$-module, and     
from (\ref{NORM}) that ${\cal N}$ is the subspace of ${\cal F}_1$
spanned by the wedge products $|\lambda\>$ labelled by $\lambda \in \DP_h - \DP$.
Therefore ${\cal F}_1/{\cal N}$ is a $A_{2n}^{(2)}$-module that can
be identified with ${\cal F}$.

In this identification one has, for 
$\lambda=(\lambda_1,\ldots,\lambda_r) \in \DP$,
\begin{equation}
P_\lambda = 2^{\sum_{i=1}^r\lfloor (\lambda_i-1)/h \rfloor} |\lambda \>.
\end{equation}
The power of $2$ comes from the fact that if $\lambda_i = kh$ for $k>0$,
and $\nu$ denotes the partition obtained from $\lambda$ by replacing $\lambda_i$
by $\nu_i = \lambda_i+1$, then it follows from (\ref{FP}), (\ref{FNP})
that $f_n P_\lambda$ contains $P_\nu$ with
coefficient 1, while 
$f_n |\lambda\>$ contains $|\nu\>$ with coefficient 2 by (\ref{ACTF}).
For later use we set
\begin{equation}\label{AHN}
a_h(\lambda) = \sum_{i=1}^r\left\lfloor {\lambda_i-1\over h} \right\rfloor \,.
\end{equation}

\section{Modular representations of $\SP_m$}

We refer the reader to \cite{B} for an up-to-date review of the
representation theory of the spin symmetric groups
and their combinatorics. 

Let $\SP_m$ be the spin symmetric group as defined by Schur \cite{S}, 
that is,
the group of order $2\,m!$ with generators $z,s_1,\ldots,s_{m-1}$
and relations $z^2=1$, $zs_i=s_iz$, $s_i^2 = z$, $(1\le i\le m-1)$,
$s_is_j=zs_js_i$ ($|i-j|\ge 2$) and $(s_is_{i+1})^3=z$
($1\le i\le m-2)$. 

On an irreducible representation of $\SP_m$, the central element $z$ has to act
by $+1$ or by $-1$. The representations for which $z=1$ are
actually linear representations of the symmetric group $\SG_m$,
and those with $z=-1$, called spin representations
correspond to two-valued representations of $\SG_m$. The irreducible spin
representations over a field of characteristic $0$
are labelled, up to association, by strict partitions
$\lambda\in \DP(m)$. More precisely, let $\DP_+(m)$ (resp. $\DP_-(m)$)
be the set of strict partitions of $m$ having an even (resp. odd) 
number of even parts. Then, to each $\lambda\in\DP_+(m)$ corresponds
a self-associate irreducible spin character $\pr{\lambda}$, and to each
$\lambda\in\DP_-(m)$ a pair of associate irreducible spin characters denoted
by $\pr{\lambda}$ and $\pr{\lambda}'$. 

According to Schur \cite{S}, the values $\pr{\lambda}(\rho)$
of the spin character $\pr{\lambda}$ on conjugacy classes of cycle-type 
$\rho=(1^{m_1},3^{m_3},\ldots )$ are given by the expansion
of the symmetric function $P_\lambda$ on the basis of power sums, 
namely
\begin{equation}
P_\lambda = \sum_\rho 2^{\lceil (\ell(\rho)-\ell(\lambda))/2\rceil}
            \pr{\lambda}(\rho) {p_\rho\over z_\rho}
\end{equation}
where $z_\rho=\prod_j j^{m_j} m_j!$ and $\ell(\lambda)$ stands for the length
of $\lambda$, that is the number of parts of $\lambda$.

For $\lambda\in\DP(m)$, one introduces the self-associate spin character
\begin{equation}
\prh{\lambda} = \cases{\pr{\lambda} & if  $\lambda\in\DP_+(m)$,\\
                     \pr{\lambda}+\pr{\lambda}'& if  $\lambda\in\DP_-(m)$.\\}
\end{equation}
The branching theorem for spin characters of Morris \cite{Mo1} implies that
if $\prh{\lambda}$ gets identified with a weight vector of ${\cal F}$
by setting
\begin{equation}\label{IDENT}
P_\lambda = 2^{\lfloor (m - \ell(\lambda))/2\rfloor} \, \prh{\lambda} ,
\end{equation}
then the $b_\infty$-operator $f = \sum_{i\ge 0} f^{\infty}_i$ implements 
the induction of self-associate spin characters
from $\SP_m$ to $\SP_{m+1}$.  
Similarly, $e= e^{\infty}_0 + 2\sum_{i>0} e^{\infty}_i$ implements
the restriction from $\SP_m$ to $\SP_{m-1}$.
Thus, the Fock space representation of $b_\infty$ may be viewed
as the sum  
${\cal F} = \bigoplus_m {\cal C}(m)$
of additive groups generated by self-associate spin characters
of $\SP_m$ in characteristic 0.
In this setting, the Chevalley generators of $b_\infty$ act as 
refined induction and restriction operators.

Now, similarly to the case $A_{n-1}^{(1)}$, the reduction from
$b_\infty$ to $A_{2n}^{(2)}$ parallels the reduction modulo $p=h=2n+1$
of representations of $\SP_m$ (from now on we assume that $h$
is an odd prime). 

More precisely, using (\ref{FP}) (\ref{FNP}) (\ref{IDENT}),
one sees immediately that the Chevalley generators $f_i$ of
$A_{2n}^{(2)}$ act on $\prh{\lambda}$ as 
the $(r,\overline{r})$-induction operators of Morris and
Yaseen $(r=n+1-i)$ \cite{MY}.
Hence the vectors of degree $m$ of 
$V(\Lambda_n) = U(A_{2n}^{(2)})^-\,|0\>$ can be identified
with linear combinations of self-associate spin characters
obtained by a sequence of $(r,\overline{r})$-inductions.
It is known from modular representation theory that the maximal
number of linearly independent self-associate projective spin characters
of $\SP_m$ in characteristic $p$ is equal to the number of partitions
of $m$ into odd summands prime to $p$.
Therefore the following result follows at once from (\ref{CHAR}).
\begin{theorem}
The  self-associate projective spin characters
of $\SP_m$ in characteristic $p$ are linear combinations 
of characters obtained by a sequence of $(r,\overline{r})$-inductions.
\end{theorem}
This was proved by Bessenrodt {\it et al.} for $p=3$ \cite{BMO}
and Andrews {\it et al.} for $p=5$ \cite{ABO}, but the question
remained open for $p\ge 7$ \cite{B}. 

Moreover, the construction of
Section~\ref{SECT4} gives an explicit basis for the space spanned
by such characters.
Denote by $\underline{A}(\mu)$ the column vector obtained
from $A(\mu)$ by reduction $q=1$ and expansion on the basis
$\prh{\lambda}$.
Then, $\underline{A}(\mu)$ is a projective character by
(\ref{defA}) and 
$\{\underline{A}(\mu)\ | \ \mu \in \DPR_p(m) \}$ is a 
basis of the $\Q$-vector space of self-associate projective spin characters
of $\SP_m$ in characteristic $p$.

These observations and the results of \cite{LLT,Ar,Gr,LT} lead us to 
formulate a conjecture relating the global basis of $V(\Lambda_n)$
and the decomposition matrices for spin characters of the groups $\SP_m$.

Let $\mu\in \DPR_p(m)$ and let $\underline{G}(\mu)$ stand for the image of the 
global basis $G(\mu)$ in ${\cal F}={\cal F}_1/{\cal N}$, that is, 
\begin{equation}
\underline{G}(\mu) = \sum_{\lambda \in \DP(m)} 
2^{b(\lambda) - a_p(\lambda)}
d_{\lambda\mu}(1) \prh{\lambda} \,,
\end{equation}
where $a_p(\lambda)$ is given by (\ref{AHN}) and 
\begin{equation}
b(\lambda)= \left\lfloor {m - \ell(\lambda)\over 2}\right\rfloor\,.
\end{equation}
Then denote by $\underline{\underline{G}}(\mu)$ the vector obtained
by factoring out the largest power of $2$ dividing the coefficients
of $\underline{G}(\mu)$ on the basis $\prh{\lambda}$.
For simplicity of notation, we shall identify $\underline{\underline{G}}(\mu)$
with the column vector of its coordinates on $\prh{\lambda}$.

Finally, let us call reduced decomposition matrix of $\SP_m$ in characteristic
$p$ the matrix obtained from the usual decomposition matrix for spin characters
by adding up pairs of associate columns and expanding the column vectors
so obtained on the basis $\prh{\lambda}$. 
This is a matrix with $|\DP(m)|$ rows and $|\DPR_p(m)|$ columns.
The definition is illustrated in Table~\ref{TAB2} and Table~\ref{TAB3}.
(Table~\ref{TAB2} is taken from \cite{MY}, except for the column labels 
which are ours and will be explained in the next section.)
\begin{table}
\caption{\label{TAB2}
The decomposition matrix of $\SP_{10}$ in characteristic 3.}
\begin{indented}
\item[]\begin{tabular}{@{}llllllll}
\br
             &(3331)&(3331)'&(4321)&(4321)'&(532)&(541)&(541)'\\
\mr
$\pr{4321}$    &0&0&1&1&0&0&0\\
$\pr{532}$     &0&0&0&0&1&0&0\\
$\pr{532}'$    &0&0&0&0&1&0&0\\
$\pr{541}$     &1&1&1&1&0&0&1\\
$\pr{541}'$    &1&1&1&1&0&1&0\\
$\pr{631}$     &2&2&1&1&0&1&1\\
$\pr{631}'$    &2&2&1&1&0&1&1\\
$\pr{64}$      &1&1&0&0&0&1&1\\
$\pr{721}$     &1&1&0&1&0&0&1\\
$\pr{721}'$    &1&1&1&0&0&1&0\\
$\pr{73}$      &1&1&1&1&0&1&1\\
$\pr{82}$      &0&0&0&0&1&0&0\\
$\pr{91}$      &1&1&1&1&0&0&0\\
$\pr{10}$      &0&1&0&0&0&0&0\\
$\pr{10}'$     &1&0&0&0&0&0&0\\
\br
\end{tabular}
\end{indented}
\end{table}

\begin{table}
\caption{\label{TAB3}
The reduced decomposition matrix of $\SP_{10}$ in characteristic 3.}
\begin{indented}
\item[]\begin{tabular}{@{}lllll}
\br
             &(3331)&(4321)&(532)&(541)\\
\mr
$\prh{4 3 2 1}$& 0 &2&0&0\\
$\prh{5 3 2}$ &0&0&1&0\\
$\prh{5 4 1}$&2&2&0&1\\
$\prh{6 3 1}$&4&2&0&2 \\
$\prh{6 4}$&2 &0&0&2\\
$\prh{7 2 1}$&2  &1 &0&1\\
$\prh{7 3}$&2 & 2&0& 2 \\
$\prh{8 2}$& 0&0&1&0\\
$\prh{9 1}$&2&2&0&0\\
$\prh{10}$&1 &0&0&0 \\
\br
\end{tabular}
\end{indented}
\end{table}

\begin{conjecture}
(i) The set of column vectors of the reduced decomposition matrix of $\SP_m$
in odd characteristic $p$ such that $p^2 > m$
coincides with $\{\underline{\underline{G}}(\mu) \ | \ \mu \in \DPR_p(m)\}$. 

(ii) For $p^2\le m$, 
the reduced decomposition matrix of $\SP_m$
is obtained by postmultiplying the matrix whose columns are 
$\underline{\underline{G}}(\mu)$ by a unitriangular matrix with 
nonnegative entries.
\end{conjecture}

Our conjecture has been checked on  the numerical tables
computed by Morris and Yaseen ($p=3$) \cite{MY} and
Yaseen ($p=5,7,11$) \cite{Ya}.
Thus, for $p=3$, $m=11$, the columns of the reduced decomposition matrix
are
\[
\underline{\underline{G}}(3332),\
\underline{\underline{G}}(4331)+\underline{\underline{G}}(641),\ 
\underline{\underline{G}}(5321),\ 
\underline{\underline{G}}(542),\ 
\underline{\underline{G}}(641).
\]

\section{Labels for irreducible modular spin characters 
and partition identities}

The labels for irreducible modular representations of symmetric
groups form a subset of the ordinary labels
\cite{JK}. It is therefore natural to look for a labelling scheme
for irreducible modular spin representations of 
$\SP_m$  using a subset of $\DP(m)$. This was accomplished for
$p=3$ by Bessenrodt {\it et al.} \cite{BMO}, who found 
that the Schur regular partitions of $m$ form a convenient
system of labels. These are the partitions 
$\lambda=(\lambda_1,\ldots,\lambda_r)$  such that
$\lambda_i-\lambda_{i+1}\ge 3$ for $i=1,\ldots,r-1$, and
$\lambda_i -\lambda_{i+1}>3$ whenever $\lambda_i\equiv 0\mod 3$.
 
In \cite{BMO}, it was also conjectured that for $p=5$, the  labels
should be  the partitions $\lambda=(\lambda_1,\ldots,\lambda_r)$
satisfying the following conditions: (1) $\lambda_i>\lambda_{i+1}$
for $i\le r-1$,
(2) $\lambda_i -\lambda_{i+2} \ge 5$ for $i\le r-2$,
(3) $\lambda_i -\lambda_{i+2} > 5$ if $\lambda_i\equiv 0\mod 5$
or if $\lambda_i+\lambda_{i+1}\equiv 0\mod 5$ for $i\le r-2$,
and (4) there are no subsequences of the following types
(for some $j\ge 0$): $(5j+3,5j+2)$, $(5j+6,5j+4,5j)$,
$(5j+5,5j+1,5j-1)$, $(5j+6,5j+5,5j,5j-1)$. 
This conjecture turned out to be equivalent to a 
$q$-series identity conjectured long ago by Andrews
in the context of extensions of the Rogers-Ramanujan identities,
and  was eventually proved by Andrews {\it et al.} \cite{ABO}.
The authors of \cite{ABO} observed however that such a labelling
scheme could not be extended to $p=7,11,13$ (see also \cite{B}).  

In terms of canonical bases, the obstruction can be understood as
follows. 
Assuming our conjecture and using the results of \cite{BMO,ABO},
one can see that for $p=3,5$, the labels of \cite{BMO} and \cite{ABO}
are exactly the partitions indexing the lowest nonzero entries
in the columns of the matrices $D_m(q) = [d_{\lambda\mu}(q)]_{\lambda,\mu\vdash m}$.
For example, in Table \ref{TAB1}, these are
$(10),(91),(82)$ and $(73)$, which are indeed the Schur regular partitions
of $10$. The problem is that for $p\ge 7$, it can happen that
two columns have the same partition indexing the lowest nonzero
entry. For example, with $p=7$ ($n=3$) and $m=21$, the two canonical basis
vectors

$ 
G(75432)
=
\ket{7 5 4 3 2}+q^{2}\ket{7 6 4 3 1}+q\ket{7 7 5 2}+q^{3}\ket{7 7 6 1}
+q^{2}\ket{8 6 4 3}+\left (q^{2}+q^{4}\right )\ket{8 6 5 2}+
q^{3}\ket{8 7 6}+q^{4}\ket{9 5 4 3}+\left (q^{4}+q^{6}\right )\ket{9 6 5 1}+
q^{5}\ket{9 7 5}
$

\noindent and

$
G(654321)
=
\ket{6 5 4 3 2 1}+q\ket{7 5 4 3 2}+ q\ket{7 6 4 3 1}+ q\ket{7 6 5 2 1}+
q^{2}\ket{7 7 4 3}+q^{2}\ket{7 7 5 2}+q^{2}\ket{7 7 6 1}+q^{3}\ket{7 7 7}+
\left (q^{3}+q^{5} \right )\ket{8 6 4 3}+
\left (q^{3}+q^{5}\right )\ket{8 6 5 2}+\left (q^{4}-q^{8}\right )\ket{8 7 6}+
\left (q^{3}+q^{5}\right )\ket{9 6 5 1}+\left (q^{4}+q^{6}\right )\ket{9 7 5}
$

\noindent have the same bottom partition $(975)$
(compare \cite{B}, end of Section 3).

On the other hand the partitions
indexing the highest nonzero entries in the columns of $D_m(q)$ are
the labels of the crystal graph (by Theorem \ref{TH}(ii)),  so that
they are necessarily distinct. Therefore, we propose to use the set 
\begin{eqnarray*}
\DPR_p(m) = &\{ \lambda = (\lambda_1,\ldots ,\lambda_r)\vdash m \ | \ 
0<\lambda_i-\lambda_{i+1}\le p \ {\rm if} \ \lambda_i \not \equiv 0 \mod p, \\
& 0 \le\lambda_i-\lambda_{i+1} < p \ {\rm if} \ \lambda_i \equiv 0 \mod p,
(1\le i \le r)\} 
\end{eqnarray*}
for labelling the irreducible spin representations of $\SP_m$
in characteristic $p$. 
Indeed its definition is equally simple for all $p$. 
Moreover, because of Theorem~\ref{TH}(iii), this labelling would be compatible with
the $p$-block structure, which can be read on the 
$\overline{p}$-cores. 
Also, it is adapted to the calculation of the vectors 
$\underline{A}(\mu)$ which give an approximation to the 
reduced decomposition matrix.

Finally, we note that since $\DPR_p$ provides the right number
of labels we have the following partition identity
\begin{equation}\label{PARTID}
\sum_{m\ge 0} | \DPR_p(m) |  t^m
=
\prod_{\scriptstyle i \ {\rm odd} \atop \scriptstyle i\not\equiv 0\mod p}
{1\over 1-t^i}
\end{equation}
which for $p=3,5$ is a counterpart to the
Schur and Andrews-Bessenrodt-Olsson identities.

This happens to be a particular case of a theorem of Andrews
and Olsson \cite{AO}.
Namely, one gets (\ref{PARTID}) by taking 
$A=\{1,2,3,\ldots ,p-1\}$ and $N=p$ in Theorem~2 of~\cite{AO}.
A combinatorial proof of a refinement of the Andrews-Olsson
partition identity has been given by Bessenrodt \cite{B1}.

One can also get a direct proof of (\ref{PARTID})
without using representation theory by 
simply considering the bijections (\ref{MAP}).

\section{Discussion}

We have used the level 1 $q$-deformed Fock spaces of Kashiwara {\it et al.}
to compute the canonical basis of the basic representation
of $U_q(A_{2n}^{(2)})$, and we have formulated a conjectural 
relation with the decomposition matrices of the spin symmetric groups 
in odd characteristic $p=2n+1$.

As in the case of $A_{n-1}^{(1)}$,
it is reasonable to expect that in general,
that is when $2n+1$ is not required to be a prime, the canonical basis 
is related to a certain family of Hecke
algebras at $(2n+1)$th roots of unity. A good candidate might be the
Hecke-Clifford superalgebra introduced by Olshanski \cite{Ol}.

The case of $2n$th roots of unity should then be related
to the Fock space representation of the affine Lie algebras
of type $D_{n+1}^{(2)}$. 
In particular we believe that the fact used by Benson \cite{Be}
and Bessenrodt-Olsson \cite{BO} that the 2-modular irreducible
characters of $\SP_m$ can be identified with the 2-modular irreducible
characters of $\SG_m$ corresponds in the realm of affine Lie algebras
to the isomorphism $D_2^{(2)} \simeq A_1^{(1)}$.

\section*{Acknowledgements}

We thank T. Miwa and A.O. Morris for stimulating discussions,
and G.E. Andrews for bringing references \cite{AO,B1} to
our attention.


\section*{References}

\begin{thebibliography}{99}

\bibitem{ABO} Andrews G E, Bessenrodt C and Olsson J B 1994 {\it
Trans. Amer. Math. Soc.} {\bf 344} 597

\bibitem{AO} Andrews G E and Olsson J B 1991 {\it J. reine angew. 
Math.} {\bf 413} 198

\bibitem{Ar} Ariki S 1996 {\it On the decomposition numbers of
the Hecke algebra of $G(m,1,n)$}, preprint

\bibitem{ANY} Ariki S, Nakajima T and Yamada H-F 1995 {\it J. Phys. A:
Math. Gen.} {\bf 28} L357
 
\bibitem{Be} Benson D 1987 {\it Proc. Sympos. Pure Math.} {\bf 47} 381

\bibitem{B1} Bessenrodt C 1991 {\it Europ. J. Combinatorics} {\bf 12} 271

\bibitem{B} Bessenrodt C 1994 {\it S\'eminaire Lotharingien de Combinatoire}
B{\bf 33}a,\\
http://cartan.u-strasbg.fr:80/\~{}slc/opapers/s33bess.html

\bibitem{BMO} Bessenrodt C, Morris A O and Olsson J B 1994 {\it
J. Algebra} {\bf 164} 146

\bibitem{BO} Bessenrodt C and Olsson J B 1993 {\it Institut fur
Experimentelle Mathematik (Essen)} preprint 17

\bibitem{DJKM1} Date E, Jimbo M, Kashiwara M and Miwa T 1982
{\it Physica D} {\bf 4} 343

\bibitem{DJKM2} Date E, Jimbo M, Kashiwara M and Miwa T 1982
{\it Publ. RIMS Kyoto Univ.} {\bf 18} 1077


\bibitem{FLOTW} Foda O, Leclerc B, Okado M, Thibon J-Y and Welsh T A
1997 {\it preprint}, q-alg/9701020
  
\bibitem{FOW} Foda O, Okado M and Warnaar S O 1996 {\it J. Math. Phys.}
{\bf 37} 965

\bibitem{G} Grojnowski I 1994 {\it
Internat. Math. Res. Notices} {\bf 5} 215

\bibitem{Gr} Grojnowski  I 1995 {\it Personal communication}

\bibitem{H} Hayashi T 1990 {\it
Commun. Math. Phys.} {\bf 127} 129

\bibitem{J} James G 1990 {\it
Proc. London Math. Soc.} {\bf 60} 225

\bibitem{JK} James G and Kerber A 1981 {\it The Representation Theory
of the Symmetric Group} (Reading Mass.: Addison-Wesley)

\bibitem{JY} Jarvis P and Yung C M 1994
{\it Lett. Math. Phys.} {\bf 30} 45

\bibitem{KKLW} Kac V G, Kazhdan D A, Lepowsky J and Wilson R L
1981 {\it Adv. Math.} {\bf 42} 83

\bibitem{KMN1} Kang S-J, Kashiwara M, Misra K C, Miwa T, Nakashima T
and Nakayashiki A 1992 {\it Int. J. Mod. Phys. A} {\bf 7} Suppl. {\bf 1A} 449

\bibitem{KMN2} Kang S-J, Kashiwara M, Misra K C, Miwa T, Nakashima T
and Nakayashiki A 1992 {\it Duke Math. J.} {\bf 68} 499

\bibitem{K} Kashiwara M 1991 {\it
Duke Math. J.} {\bf 63} 465

\bibitem{KMPY} Kashiwara M, Miwa T, Petersen J-U H and Yung C M 1996
{\it Selecta Math.} {\bf 2} 415 

\bibitem{KMS} Kashiwara M, Miwa T and Stern E 1996, {\it
Selecta Math.} {\bf 1}  787

\bibitem{LLT} Lascoux A, Leclerc B and Thibon J-Y 1996 {\it
Commun. Math. Phys.} {\bf 181} 205

\bibitem{LLTrib} Lascoux A, Leclerc B and Thibon J-Y 1997 {\it
J. Math. Phys.} {\bf 38} (2) 1041 

\bibitem{LT} Leclerc B and Thibon J-Y 1996 {\it
Internat. Math. Res. Notices} {\bf 9} 447

\bibitem{MM} Misra K C and Miwa T 1990 {\it
Commun. Math. Phys.} {\bf 134} 79

\bibitem{Mo1} Morris A O 1965 {\it
Can. J. Math.}  {\bf 17} 543

\bibitem{MY1} Morris A O and Yaseen A K 1986 {\it
Math. Proc. Camb. Phil. Soc.} {\bf 99} 23 

\bibitem{MY} Morris A O and Yaseen A K 1988 {\it
Proc. Royal Soc. Edinburgh} {\bf 108A} 145

\bibitem{NY} Nakajima T and Yamada H-F 1994 {\it J. Phys. A: Math. Gen.} {\bf 27} L171

\bibitem{Ol} Olshanski G I 1992 {\it Lett. Math. Phys.} {\bf 24} 93

\bibitem{S} Schur I 1911 {\it J. Reine Ang. Math.} {\bf 139} 155

\bibitem{Soe1} Soergel W 1997 {\it Represent. Theory} {\bf 1} 83

\bibitem{Soe2} Soergel W 1997 {\it Represent. Theory} {\bf 1} 115

\bibitem{St} Stern E 1995 {\it
Internat. Math. Res. Notices} {\bf 4} 201

\bibitem{Ya} Yaseen A K 1987 {\it Thesis}, University of Wales

\bibitem{You} You Y 1989 {\it Adv. Ser. in Math. Phys. }
{\bf 7} 449

\end{thebibliography}
\end{document}